\documentclass[12pt]{article}
\usepackage{graphicx}
 \usepackage{mathptmx}      
\usepackage{latexsym,amsmath,amssymb,amsfonts,amsthm}
\newcommand{\oM}{\overline{M}}

\newcommand{\sm}[1]{\mbox{\small $#1$}}

\newcommand{\La}[1]{\mbox{\Large $#1$}}

\newcommand{\ha}{\sm{\frac{1}{2}}}

\newcommand{\h}{\mbox{\lie h}}                  

\newfont{\lie}{eufm10 at 12pt}
\newfont{\field}{msbm10 at 11pt}
\newtheorem{theorem}{Theorem}
\newtheorem{lemma}{Lemma}
\newtheorem{corollary}{Corollary}
\newtheorem{proposition}{Proposition}
\newtheorem{definition}{Definition}
\begin{document}

\title{\normalsize  \bf   ON THE SPECTRUM  OF
$\bar{X}$-BOUNDED  MINIMAL
SUBMANIFOLDS }

\markright{\sl\hfill   I.\ Salavessa \hfill}

\author{Isabel M.C.\ Salavessa
}
\date{}
\protect\footnotetext{\!\!\!\!\!\!\!\!\!\!\!\!\! {\bf MSC 2000:}
Primary: 53C40; 58C40
 Partially supported by FCT through the Plurianual of CFIF.}
\maketitle ~~~\\[-15mm]

\begin{quotation}\noindent
{
\footnotesize  Centro de F\'{\i}sica das Interac\c{c}\~{o}es
Fundamentais, Instituto Superior T\'{e}cnico, Technical University
of Lisbon, Edif\'{\i}cio Ci\^{e}ncia, Piso 3, Av.\ Rovisco Pais,
P-1049-001 Lisboa, Portugal;~~
e-mail: isabel.salavessa@ist.utl.pt}\\[5mm]
\end{quotation}
\begin{abstract}
\noindent
We prove, under a certain  boundedness condition at infinity
of a  $(\bar{X}^{\top},\bar{X}^{\bot})$ component of
the second fundamental form, 
the vanishing of the essential spectrum
of  a complete minimal $\bar{X}$-bounded and $\bar{X}$-properly immersed
submanifold on a Riemannian manifold endowed with a strongly convex
 vector field $\bar{X}$.  The same conclusion also holds for
 any complete minimal $h$-bounded and $h$-properly
immersed submanifold  that lies in a open set of a Riemannian manifold
$\oM$ supporting a nonnegative strictly convex function  $h$.
This  extends a recent result
of Bessa, Jorge and Montenegro on  the spectrum of Martin-Morales minimal
surfaces. Our proof  uses as  main tool an extension
of Barta's theorem given in \cite{BM}.
\end{abstract}

\section{Introduction and main results}
\noindent
Since Calabi in 1965 \cite{Ca} conjectured that complete minimal hypersurfaces
in Euclidean spaces are unbounded, some answers have been given,
with a positive answer by Colding and Minicozzi \cite{CM}
  for the case of embedded surfaces, 
and a negative answer with the counterexamples given by  
Nadirashvili \cite{Na}  and by  Martin and Morales \cite{MM1,MM2}
for the case of immersed nonembedded surfaces.
This conjecture also motivates many other related problems in
more general ambient spaces, for instance, 
on the topological and geometrical properties of minimal
submanifolds that are   bounded or not, or on the search of  conditions
for a submanifold to be unbounded. 
In \cite{BJM} the structure of
the spectrum of the Martin-Morales surfaces is studied, 
namely it is proved that
 complete bounded  minimal properly immersed submanifolds of
the unit open ball of $\mathbb{R}^n$ 
must have pure point spectrum. 

In this note we extend the above result of Bessa, Jorge and Montenegro
to an ambient space carrying an almost conformal vector field
$\bar{X}$, a concept introduced in
(\cite{RS,S0}). On a Riemannian $(m+n)$-dimensional
 manifold $(\bar{M}, \bar{g})$
 we say a vector field  $\bar{X}$ is almost conformal  if
\begin{equation}
 2\alpha \bar{g}\leq  L_{\bar{X}} \bar{g} \leq  2\beta \bar{g}
\end{equation}
where $+\infty\geq \beta\geq \alpha > 0$ are constants, and
$L_{\bar{X}}\bar{g}(\bar{Y},\bar{Z})=\bar{g}(\bar{\nabla}_{\bar{Y}}\bar{X},
\bar{Z })+\bar{g}(\bar{\nabla}_{\bar{Z}}\bar{X},
\bar{Y })$, where $\bar{\nabla}$ is the Levi-Civita connection
of $(\oM, \bar{g})$. 
If we allow  $\beta=+\infty$, in this case
$\bar{X}$ is named  by strongly convex.

An example of almost conformal vector field in a complete Riemannian manifold
$\oM$ is  the position vector field
$\ha \bar{\nabla} r^2 = r \frac{\partial }{\partial r}$ on
a geodesic ball of $\oM$ of radius $R$ and center $\bar{p}$ 
that does not intercept the cut locus at $\bar{p}$ and
$\sqrt{\kappa^+}R<\pi/2$ with
$\kappa^+=\kappa^+(R)=\max\{0,\sup_{B_R(\bar{p})}\bar{K}\}$,  
where $\bar{K}$ are the sectional curvatures ̄
of $\bar{M}$ and $r$ is the distance function on $\oM$ to a given point. 
In this case $\alpha$ and $\beta$ are well defined functions 
$\alpha=\alpha_{\kappa^+}(R)$, $\beta=\alpha_{\kappa^-}(R)$,
of $R$, $\kappa^+$, and $\kappa^-=\kappa^-(R)=
\min\{0,\inf_{B_R(\bar{p})}\bar{K}\}$
where
\begin{equation}
\alpha_{\kappa}(R)=\left\{\begin{array}{lll}
R\sqrt{\kappa}\cot(\sqrt{\kappa}R)&\mbox{for}~ 0\leq R< \pi/2\sqrt{\kappa},
~\mbox{when}~ \kappa>0\\
1& \mbox{for}~0\leq R< +\infty,~
\mbox{when}~ \kappa=0\\
R\sqrt{-\kappa}{\coth}(\sqrt{-\kappa}R)&\mbox{for}~ 0\leq R<
+\infty, ~
\mbox{when}~ \kappa <0. \end{array}\right.
\end{equation}
A strictly convex function f on $\oM$
 with $\rm{Hess} ~f \geq  \alpha \bar{g}$
 defines a strongly convex vector field $\bar{\nabla} f$. 
Positive homothetic non-Killing vector fields are almost conformal.
In $\mathbb{R}^{m+n}$ the position vector field $\bar{X}_x = x$
 is such an example. A particular feature of strongly convex
 vector fields, is that the norm $\|\bar{X}\|$ must take its maximum
on the boundary of compact domains (see proposition 1).
Therefore $\bar{X}$ cannot be globally defined
on a compact manifold $\oM$ without boundary.

Strongly convex vector fields have a role on isoperimetric
inequalities for an immersed $m$-dimensional 
submanifold $F:M\to \oM$, $m\geq 2$, involving
the the mean curvature $H$. The Cheeger constant of $M$
is defined by $\h(M)=\inf_D{A(\partial D)}/{V(D)}$, where $D$ runs
 over all compact domains $D$ of $\oM$ with picewise smooth boundary
$\partial D\subset M$ of respective volume $V(D)$ and area $A(\partial D)$.
We recall the following inequality \cite{LS}:
\begin{equation}
 (\sup_M \|\bar{X}_F\|)^{-1}\leq \frac{1}{\alpha}
\left(\frac{1}{m}\h(M) +\sup_M\|H\|\right)
\end{equation}
where $\bar{X}_F$ denotes $\bar{X}$ along $F$. 
Let $\bar{X}^{\top}$ and $\bar{X}^{\bot}$ denote
the orthogonal projections of $\bar{X}_F$ onto $TM$, and
the normal bundle $NM$ respectively.
We remark that,
 following the proof in \cite{LS} we see that if $F$ is minimal
we have a sharper inequality:
\begin{equation}
 (\sup_M \|\bar{X}^{\top}\|)^{-1}
\leq \frac{1}{\alpha}\frac{1}{m}\h(M).
\end{equation}
We note that $\bar{X}_F$  ($\bar{X}^{\top}$ resp.)
cannot vanish identically for any (minimal resp.) immersion $F$ 
 (see lemmas 1 and 2).
In the case $\oM$ is the Euclidean space with the position vector
field, $\|\bar{X}^{\top}\|\leq \|\bar{X}_F\|=\|F\|$.
This leads to the following conclusion:
\begin{theorem}[\cite{LS}] 
If $\bar{X}$ is a strongly convex vector field on a neighbourhood
of a  minimal submanifold $F:M\to \bar{M}$ with zero Cheeger
constant, then  $\bar{X}^{\top}$ is unbounded. In the particular case
$\oM=\mathbb{R}^{n+m}$, $F$ is unbounded.
\end{theorem}

We  recall the following inequality due to Cheeger \cite{Che},
$$ \h^2(D)\leq 4\lambda(D)$$
where $\lambda(D)$ is the fundamental tone of a normal domain $D$
 in $M$. For normal bounded domains, $\lambda(D)$ is the first eigenvalue
for the boundary Dirichlet problem.
The  Rayleigh characterization of the fundamental tone of 
any open  domain $D$ of $M$ is given by
$$\lambda(D)=\inf\left\{ \frac{\int_D\|\nabla f\|^2}{\int_Df^2} 
: f\in L^2_{1,0}(D)\right\}$$
where $L^2_{1,0}(D)$ is the completation of $C^{\infty}_0(D)$
for the norm $\|\phi\|^2=\int_M\phi^2+ \|\nabla\phi\|^2$.
Thus, 
if $M$ is complete noncompact, 
$\lambda(M)=\lim_R\lambda(D_R)$ 
and $\h(M)\leq \h(D_R)$, where
$D_R$  is an exhaustion  sequence of bounded 
domains of $M$ with smooth boundary in $M$. 
Therefore, from the above inequalities we have the following estimate
for $M$ a bounded domain (possibly with  boundary) or a complete Riemannian
manifold
\begin{eqnarray}
 (\sup_M \|\bar{X}_F\|)^{-1}&\leq &\frac{1}{\alpha}
\La{(}\frac{2}{m}\sqrt{\lambda(M)} +\sup_M\|H\|\La{)} \nonumber\\
(\sup_M \|\bar{X}^{\top}\|)^{-1} &\leq&  \frac{1}{\alpha}
\frac{2}{m}\sqrt{\lambda(M)}, ~~~~\mbox{if}~M~\mbox{is~minimal}.
\end{eqnarray}
\begin{definition} Given a vector field $\bar{X}$ of
$\oM$, an immersed submanifold $F:M\to \oM$ 
is said  $\bar{X}$-bounded  if $ \sup_M\|\bar{X}^{\top}\|<+\infty$.
If  $\sup_M\|\bar{X}^{\top}\|$ is not achieved,
then $F$ is said
$\bar{X}$-proper,
if  $\|\bar{X}^{\top}\|:M\to [0, \sup_M\|\bar{X}^{\top}\|)$ is
a proper map. 
\end{definition}
\noindent
We will see in proposition 1 that if $M$ is minimal 
and $\bar{X}$
is strongly convex, then $\sup_M\|\bar{X}^{\top}\|$ is not achieved
(in $M$) if  condition (6) below holds.
\noindent
Note that if $\bar{X}$ is the position vector field of $\oM$,
$\|\bar{X}_{p}^{\top}\|\leq \|\bar{X}_{F(p)}\|=r(F(p))$. This implies
$\bar{X}$-boundedness is a weaker concept then the usual
boundedness of   $M$  in $\oM$.
For example, the spiral curve in $\mathbb{R}^2$,
$\gamma(t)=a e^{tb}(\cos (e^{abt}), \sin(e^{abt}))$  with $a>1$ and $ b>0$
 constants, is $\bar{X}$-bounded but unbounded in the
usual sense.
On the other hand
$\bar{X}$-properness might be a stronger concept than the usual
properness of an immersion. We also remark that
if $\bar{X}^{\bot}=0$ along all $M$,
then  $\bar{\nabla} r$ restricted to $M$ is a vector field on $M$.
If  $r$ is the distance function on $\oM$ from a fixed point
 $p\in M$, we  see that
 (unit) geodesics of $\oM$ starting at $p$ ( that are the integral
curves of $\bar{\nabla} r$) lie in $M$. In this case $n=0$.\\

Next we state our main theorems:
\begin{theorem} Let $F:M\to \bar{M}$ be a  complete minimal immersion
that is  $\bar{X}$-bounded with $\sup_M \|\bar{X}^{\top}\|= R$, 
where $\bar{X}$ is a strongly convex vector field of $\oM$ defined on a 
neighbourhood of $M$, then:\\[3mm]
$(1)$~~$2\sqrt{\lambda(M)}\geq \h(M)\geq \frac{m \alpha }{R}$.\\[2mm]
$(2)$ Furthermore, if  the second fundamental form
$B$ of $M$ satisfies  at points $p\in M$ with $\|\bar{X}^{\top}\|$ 
sufficiently close to $R$,
\begin{equation}
|\bar{g}(B(\bar{X}^{\top},\bar{X}^{\top}),\bar{X}^{\bot})|\leq  
\alpha'
\|\bar{X}^{\top}\|^2, 
\end{equation}
for some nonnegative constant $\alpha'<\alpha$, 
and if  $M$ is $\bar{X}$-proper,
 then
the spectrum of $M$ is a pure point spectrum.
\end{theorem}
\noindent
The condition (6) does not mean $\|B\|$ is bounded, even in the case
$\bar{X}$ is the position vector field $r\frac{\partial}{\partial r}$.
In theorem 5 (section 2) we will see that boundedness of the second 
fundamental form is, in general,  not a compatible condition
 with the boundedness of a complete minimal
submanifold, 
for ambient spaces with sectional curvature bounded
from below.
Moreover, 
for the particular case  of $\bar{X}$ being the gradient of
a nonnegative convex smooth function $h:\oM\to [0, +\infty)$ we can remove
the  boundedness condition  (6) of theorem 2, 
if we  adapt our definition of boundedness and of properness:
$F$ is $h$-bounded if  $\sup_M h\circ F=R<+\infty$, and 
is
$h$-proper if $\sup_M h\circ F$ is not achieved and 
$h\circ F:M\to [0, R)$ is a proper function.
We also will see in proposition 1 that $\sup_M h\circ F$ cannot be achieved
for $F$ a  minimal immersion.
\begin{theorem} Let $h:\oM\to [0, +\infty)$ be a nonnegative convex 
smooth function and 
$F:M\to \bar{M}$  a complete minimal
immersion that is  $h$-bounded.
If $F$ is $h$-proper,
 then the spectrum of $M$ is
a pure point spectrum.
\end{theorem}
\noindent
The above case contains the next example, when $h= \ha r^2$, where
$r$ is the
distance function to a point $\bar{p}$ in $\oM$. Note that if $F$ is
$h$-bounded, then it is also $\bar{X}$-bounded, for $X$ the position
vector field, and the concept of $h$-bounded ($h$-porper resp.)
is equivalent to usual boundedness (properness resp.). Next corollary
is a corollary of theorem 2 (1) and theorem 3:
\begin{corollary}  If ~$F:M\to\oM$ is a complete bounded minimal submanifold
with $F(M)$ lying in a  open geodesic ball $ B_R(\bar{p})$ of
$\oM$, and  $R$ is in the conditions given in $(2)$,
then $2\sqrt{\lambda(M)}\geq \h(M)\geq \frac{m\alpha}{R}$,
where $\alpha=\alpha_{k^+}(R)$.
Furthermore, if $F$ is a proper immersion into $B_R(\bar{p})$,
 then the spectrum of  $M$ is a pure point spectrum.
\end{corollary}
\begin{corollary} If ~$F:M\to \oM$ is a complete bounded minimal submanifold
properly immersed in $B_R(\bar{p})$, and  $\oM$ is a complete
Riemannian manifold with $\bar{K}\leq 0$, then $2\sqrt{\lambda(M)}\geq
\h(M)\geq \frac{m}{R}$ and the spectrum of $M$ is a pure point spectrum.
\end{corollary}
\noindent
The later corollaries are  straightforward  generalizations of  \cite{BJM}.
 Donnelly in \cite{Donn} 
 proved the existence of a non-empty essential spectrum for  negatively
curved manifolds under certain conditions.  This result and corollary 2
gives next corollary:
\begin{corollary} 
There is no complete simply connected 
minimal surface 
 $F:M^2\to \oM$ 
properly immersed into a geodesic ball $B_R(\bar{p})$ of a space form
$\oM$ of constant sectional curvature $ \bar{K}<0$, 
 and
satisfying $\|B\|^2\to c$ at infinity, for any nonnegative
finite constant $c$.
\end{corollary}
As we have announced above, in theorem 5  we will see  this conclusion can be 
extended to a considerably more general setting, where we do not need
to use spectral theory to prove it, but  a generalized Liouville-type 
 result due to  Ranjbar-Motlagh \cite{R-M}.
\\

\noindent
An application of a hessian  comparison theorem
for the distance function to a totally convex submanifold due to Kasue 
\cite{Ka}  give us the following theorem:
\begin{theorem}  Let  $\oM$ be a connected complete Riemannian manifold with
nonnegative sectional curvature and $\Sigma$  a totally convex submanifold
of dimension $d\geq n$ 
that is a closed subset of $\oM$, and  let $h=\ha \rho^2$,
where $\rho$ is the distance function in $\oM$ to $\Sigma$. If
$F:M\to \oM$ is a complete minimal immersed submanifold such that for any
$p\in M\backslash F^{-1}(\Sigma)$, $\|(\sigma'_{F(p)}(l))^{\top}\|^2
\geq {\alpha}$,
where $0<\alpha\leq 1$ is a constant and
 $\sigma_{F(p)}:[0,l]\to \oM$ is the unique  geodesic normal to
$\Sigma$ that satisfies $\sigma_{F(p)}(0)\in \Sigma$ and $\sigma_{F(p)}(l)=
F(p)$, then:\\[2mm]
$(1)$ ~$2\sqrt{\lambda(M)}\geq \h(M)\geq \frac{m\alpha}{\sup_M \rho\circ F}$.
In particular, if $M$ has zero Cheeger constant, then $\rho\circ F$ is
unbounded.\\[1mm]
$(2)$ If $M$ is $h$-bounded and $h$-properly immersed,
then  $M$ has pure point spectrum only.
\end{theorem}
\noindent
In the last section we apply this general result to submanifolds
of a product of Riemannian manifolds.

\section{Some inequalities for minimal submanifolds}
Let $\bar{X}$ be an almost conformal vector field of $\oM$, and
$F:M\to \oM$ an immersion of a $m$-dimensional submanifold with
second fundamental form $B:\bigodot^2 TM\to NM$, where
$NM$ is the normal bundle of $M$. 
We give to $M$ the induced Riemannian metric $g=F^*\bar{g}$
and the corresponding Levi Civita connection $\nabla$. 
We denote by $(\cdot)^\top$ and $(\cdot)^{\bot}$ the orthogonal
projections of $T_{F(p)}\oM$ onto $T_pM\equiv dF_p(T_pM)$ and $NM_p$
respectively.
We have for $X,Y$ vector fields on $M$, 
$ \nabla_X Y= (\bar{\nabla}_X Y)^{\top}$ and $B(X,Y)=(\bar{\nabla}_X Y)^{\bot}.
$
The mean curvature of $M$ is the normal vector
given by $H=\frac{1}{m}\rm{trace}_gB$.
The projection $\bar{X}^{\top}$ defines a vector field on $M$, and
$\bar{X}^{\bot}$ a section of the normal bundle. Since $\bar{X}_F=
\bar{X}^{\top}+\bar{X}^{\bot}$, 
an elementary computation
gives 
\begin{lemma} For $Y, Z\in T_pM$, ~$
L_{\bar{X}^{\top}}{g}(Y,Z) = L_{\bar{X}}\bar{g}(Y,Z)
+2\bar{g}(B(Y,Z),\bar{X}^{\bot}).$ In particular, $\bar{X}_F$ cannot 
vanish everywhere in any open domain of $M$.
\end{lemma}
\begin{lemma}
$(1)$ ~$m\alpha  +m\bar{g}(H, \bar{X}^{\bot})\leq 
{\rm div}_g(\bar{X}^{\top})\leq m\beta  +m\bar{g}(H, \bar{X}^{\bot})$.
If ~$F$ is minimal then
~$ m\alpha \leq {\rm div}_g(\bar{X}^{\top})\leq m\beta$, and
$\bar{X}^{\top}$  cannot 
vanish everywhere in any open domain of $M$.\\[2mm]
$(2)$ 
~$g(\nabla \|\bar{X}^{\top}\|,\bar{X}^{\top})\geq 
\alpha \|\bar{X}^{\top}\|
+\frac{1}{\|\bar{X}^{\top}\|}
\bar{g}(B(\bar{X}^{\top},\bar{X}^{\top}),\bar{X}^{\bot})$.
\end{lemma}
\begin{proof} 
 Let
$e_i$ be an o.n. basis of $T_pM$. At $p$, 
 ${\rm div}_g(\bar{X}^{\top})
= \sum_i\ha L_{\bar{X}^{\top}}{g}(e_i,e_i)$,
and an application of previous lemma gives (1) as well (2) since
$$
g(\nabla \|\bar{X}^{\top}\|,\bar{X}^{\top})=\sum_i
\frac{1}{ \|\bar{X}^{\top}
\|}g(\nabla_{e_i}\bar{X}^{\top}, \bar{X}^{\top})
g(e_i, \bar{X}^{\top})
=\frac{1}{2\|\bar{X}^{\top}\|}
 L_{\bar{X}^{\top}}g(\bar{X}^{\top},\bar{X}^{\top}).
$$
\end{proof}
\begin{proposition} If $\bar{X}$ is strongly convex, then:\\[2mm]
$(1)$~ For any  bounded domain $D$ of $\oM$ the norm $\|\bar{X}\|$
takes its maximum on the boundary $\partial D$.\\[1mm]
$(2)$~ If $F$ is a minimal immersion and (6) holds,
then the supremum of $\|\bar{X}^{\top}\|$
cannot be achieved. In particular $M$ cannot be compact without
boundary (closed).\\[1mm]
$(3)$~If $\bar{X}=\nabla h~$ for a smooth nonnegative convex function $h:\oM\to
\mathbb{R}$ and $F:M\to \oM$ is a minimal
submanifold, then the supremum of $h\circ F$ cannot be achieved.
In particular $M$ cannot be closed.
\end{proposition}
\begin{proof}
From the inequality $\bar{g}(\bar{\nabla}\|\bar{X}\|^2, \bar{X})=2\bar{g}
(\bar{\nabla}_{\bar{X}}\bar{X},\bar{X})
\geq 2\alpha\|\bar{X}\|^2$, all critical points of $\|\bar{X}\|^2$ are
vanishing points. This proves (1). To prove (2) we assume  a maximum 
point $p_0$ of $\|\bar{X}^{\top}\|$ exists. Then at $p_0$
 we may take $e_1=\bar{X}^{\top}/\|\bar{X}^{\top}\|$, and we have
by lemma 1 and (6)
\begin{eqnarray*}
0&=&\|\nabla\|\bar{X}^{\top}\|^2\|^2= 4\sm{\sum}_i|g(\bar{\nabla}_{e_i}
\bar{X}^{\top},\bar{X}^{\top})|^2\geq 4|g(\bar{\nabla}_{\bar{X}^{\top}}
\bar{X}^{\top}, \bar{X}^{\top})|^2 \|\bar{X}^{\top}\|^{-2}\\
&=& (L_{\bar{X}^{\top}}g(\bar{X}^{\top},\bar{X}^{\top}))^2
\|\bar{X}^{\top}\|^{-2}
\geq (2\alpha \|\bar{X}^{\top}\|^2
+2\bar{g}(B(\bar{X}^{\top},\bar{X}^{\top}),\bar{X}^{\bot}))^2
\|\bar{X}^{\top}\|^{-2}\\
&\geq& C^2
\end{eqnarray*}
where $C=2(\alpha-\alpha')$, 
what is impossible.
Finally we prove (3). A maximum point $p_0$ of $h\circ F$ satisfies
$\Delta (h\circ F)(p_0)\leq 0$, what contradicts
\begin{equation}
\Delta (h\circ F)_{p}=\sum_i ({\rm Hess} \, h)_{F(p)}(dF(e_i), dF(e_i))
+m\bar{g}(\bar{\nabla}h_{F(p)},H)
\geq m\alpha.
\end{equation}
\end{proof}
 \begin{theorem}
If $\overline{K}$ is bounded from  bellow, and $F:M\to \oM$ is
any complete immersed minimal submanifold with  bounded
second fundamental form, then for any nonnegative
strictly convex function $h:\oM\to [0,+\infty)$ defined in a
neighbourhood of $F(M)$, $F$ is $h$-unbounded.
\end{theorem}
\begin{proof}
Let us assume there exists a complete $h$-bounded 
immersion with $\|B\|^2\leq b$, $b$ a nonnegative constant. 
By Gauss equation, 
the Ricci tensor of $M$ is  bounded from below. Indeed, if $Y\in T_pM$
is a unit vector,
\begin{eqnarray*}
Ricci(Y,Y) &=& {\sum}_i\bar{K}(Y,e_i)
+ m\bar{g}(H,B(Y,Y))-\bar{g}(B(Y,e_i),B(Y,e_i))\\
&\geq&
m\,( inf_{\oM}\bar{K})-mb- b.
\end{eqnarray*}
Furthermore,
(7) holds for $F$ minimal imersion.
Then  theorem 2.1 of \cite{R-M} gives us 
$\lim\sup_{r_M(p)\to+\infty} \frac
{h\circ F(p)}{r_M(p)}\geq  C$, where $C$ is a positive constant that
depends on $m, b, \alpha$ and a lower bound of $\bar{K}$. 
This contradicts the assumption of $h\circ F$ to be bounded.
\end{proof}
\noindent
\\
Bessa and Montenegro defined in \cite{BM} a  quantity on a domain $D$
(bounded or not) of $M$, that here  we denote by $c(D)$
$$c(D) =\sup_X \left( \inf_D ({\rm div}_g X -\|X\|^2) \right)$$
where $X$ runs over all vector fields on $D$ locally integrable and
with a weak divergence. We denote by $c(X)={\rm div}_g X-\|X\|^2$.
\begin{proposition}[\cite{BM}] 
$\lambda(D)\geq c(D)$, 
with equality if $D$ has compact closure with smooth boundary.
\end{proposition}
\noindent
 Assume $\sup_M\|\bar{X}^{\top}\|=R<+\infty$ and (6) holds.
Set $C=2(\alpha-\alpha')$.
For each  $\epsilon>0$  sufficiently small constant we consider the domain 
$$D_{\epsilon}=\{p\in M:R^2> \|\bar{X}^{\top}\|^2>R^2-\epsilon^2\}.$$
\begin{proposition} If $F$ is a minimal submanifold  
and  (6) holds, then for any $0<\epsilon<R$ sufficiently small,
$$\lambda(D_{\epsilon})\geq \frac{mC\alpha}{\epsilon^2}.  $$
\end{proposition}
\noindent
\begin{proof} 
We define the function $f:[\sqrt{R^2-\epsilon^2}, R)\to [\epsilon,+\infty)$, 
$f(s)=\frac{C}{R^2-s^2}$, and  the smooth
vector field on $D_{\epsilon}$, $X=f(t)\bar{X}^{\top}$,
where $t= \|\bar{X}^{\top}\|$.
Using   lemma 2, we have
\begin{eqnarray*}
c(X) &=& f(t)div_g(\bar{X}^{\top})
+g(\nabla (f(t)), \bar{X}^{\top})
-f(t)^2t^2\\
&\geq & f(t)m\alpha + f'(t)(\alpha t +\frac{1}{t}g(B(\bar{X}^{\top},
\bar{X}^{\top}),\bar{X}^{\bot}))- f^2(t) t^2.
\end{eqnarray*}
Note that $f'(s)$ and $f^2(s)$ go faster to $+\infty$  then $f(s)$, when
$s \to R$. Then we have to require
$f'(t)(\alpha t +\frac{1}{t}g(B(\bar{X}^{\top},\bar{X}^{\top}),
\bar{X}^{\bot}))- f^2(t) t^2\geq 0,$
that holds under condition (6). In this case,
$$c(X)\geq \frac{C m\alpha}{R^2-t^2}\geq \frac{C m\alpha}{\epsilon^2}.$$
Now proposition 2  gives the lower bound for $\lambda(D_{\epsilon})$.
\end{proof}
\section{Proof of theorems 2 and 3}
Let $M$ be a complete noncompact  $m$-dimensional Riemannian manifold, with
Laplacian operator $\Delta$ acting on the  domain $\mathcal{D}$ of  $L^2(M)$,
where $\Delta \phi \in L^2$ for any $\phi\in \mathcal{D}$. 
The spectrum of $-\Delta$ 
decomposes
as $\sigma(M)=\sigma_{p}(M)\cup \sigma_{ess}(M)\subset 
[\lambda(M),\infty)$, where 
$\sigma_{p}(M)$ is the pure point spectrum of isolated finite multiplicity
eigenvalues, and $\sigma_{ess}(M)$ is the essential spectrum.
The decomposition principle of \cite{DL} states that 
$M$ and $M\backslash K$ have the same essential spectrum, as long as
$K$ is a compact domain of $M$ with boundary. \\[4mm]

\noindent
\em Proof of theorem 2. \em 
 (1) is immediate from (4) and the Cheeger inequality. (2)
We can take a sequence $\epsilon_k\to 0$ such that $\sqrt{R^2-{\epsilon_k}^2}$ 
are regular values of $\|\bar{X}_F^{\top}\|$. Since  $F$ is $\bar{X}$-proper,
the sets $K_{\epsilon_k}=M\backslash D_{\epsilon_k}$ 
are compact with smooth boundary.
 As in \cite{BJM} we prove the theorem
by showing that $\lambda(D_{\epsilon_k})\to +\infty$ when $k\to +\infty$,
what proves that $\sigma_{ess}(M)=\emptyset$. This is the
case by proposition 3. \qed 
\vspace*{5mm}
\noindent
\em Proof of theorem 3. \em 
In this case we take  the domain of $M$,
$D_{\epsilon}=\{p\in M: R>h\circ F>R-\epsilon\}$, 
and the vector field defined  on $D_{\epsilon}$ given by
$X={\nabla(h\circ F)}/{(R- h\circ F)}$. 
Then using (7), 
$$
c(X)_p = \frac{\Delta(h\circ F)(p)}
{(R- h( F(p)))}=\frac{\sum_i(\rm{Hess} \, h)_{F(p)}(dF(e_i),dF(e_i))}
{(R- h( F(p)))}\geq \frac{m \alpha}{\epsilon},
$$
and so $\lambda(D_{\epsilon})\to +\infty$ when $\epsilon\to 0$.\qed
\\[4mm]
\noindent
\em Proof of corollary 3  \em   ~Assume such immersion exists
with  $\|B\|^2\to c$ at infinity, $c\geq 0$ a finite constant. By the Gauss 
equation the sectional curvature
of $M$ satisfy $ K=\bar{K}-\|B\|^2$. Then $M$ has negative
sectional curvature and  $K\to \bar{K}-c <0$
at infinity. By a result of Donnelly \cite{Donn} the essential spectrum of $M$  
consists of the half line $[(-\bar{K}+c)/4,+\infty)$ contradicting corollary 2.
\qed
\section{Ambient space with a  totally convex set} 
\begin{definition} $(1)$~We say a vector field 
$\bar{X}$ of $\oM$ is  almost trace-conformal (strongly trace-convex
resp.) along $M$ if $2m\alpha  \leq {\rm Trace}_g~F^*L_{\bar{X}}\bar{g}\leq 
2\beta m$ (with $\beta=+\infty$ resp.), where $\beta\geq \alpha>0$
are constants.\\[1mm]
$(2)$~ We say that a function $h:\oM\to [0,+\infty)$ is
 strictly trace-convex along $M$ if for some positive constant
$\alpha$,
${\rm Trace}_g F^*({\rm Hess}~h)\geq m\alpha$.
\end{definition}
\noindent
It is elementary to verify next theorem, following the previous proofs:
\begin{theorem}
In the  weaker conditions of  definitions 2 and  1, the inequality
 $(4)$ still holds as well the conclusions in theorems 1, 2 and 3.
\end{theorem}
\noindent
A subset  $\Sigma$ of 
$\oM$  is said to be totally convex if
it contains any geodesic connecting two points of $\Sigma$.
If $\Sigma$ is a submanifold  that is a closed subset of
$\oM$,
the hessian of the function $h= \ha \rho^2$, 
where $\rho$ is the distance function
in $\oM$ to $\Sigma$,  satisfies the following  comparison theorem:
\begin{theorem}[\cite{Ka}] If ~$\oM$ is a connected complete Riemannian manifold
with nonpositive sectional curvature and $\Sigma$ is a totally convex
submanifold of dimension $d$ that is a closed subset of $\oM$, 
then for any $Y\in T_q\oM$, 
 $q\notin \Sigma$,
$$({\rm Hess}~ h)_q(Y,Y)\geq \bar{g}(\sigma'_q(l),Y)^2$$
where $\sigma_q:[0,l]\to \oM$ is the unique unit  geodesic normal to $\Sigma$
that satisfies  $\sigma_q(0)\in \Sigma$ and $\sigma_q(l)=q$. 
\end{theorem}
\noindent 
In \cite{Ka} the condition on $\sigma$ is that it must satisfy
$t=\rho(\sigma(t))$, but this is equivalent to
$\sigma$ meets $\Sigma$ orthogonally (see \cite{Gray} chapter 2).
\begin{proposition} If  $\oM$ is in the conditions of theorem 7 and 
$F:M\to \oM$ is a complete minimal immersed submanifold such that for any
$p\in M\backslash F^{-1}(\Sigma)$, 
$\|(\sigma'_{F(p)}(l))^{\top}\|^2\geq{\alpha}$,
where $0<\alpha\leq 1$ is a constant (in particular  $d\geq n$), 
then $h$ 
is strictly trace-convex
along $M$, $\sup_M  \rho\circ F$ cannot be achieved, and
$$(\sup_M ~\rho\circ F)^{-1}\leq \frac{1}{m}\frac{1}{\alpha}\h(M).$$
 \end{proposition}
\begin{proof}
From theorem 7, 
$$\sum_i({\rm Hess}~h)_{F(p)}(dF(e_i),dF(e_i))\geq 
\sum_i(\bar{g}(\sigma'_{F(p)}(l),dF(e_i)))^2
=\|(\sigma'_{F(p)}(l))^{\top}\|^2,
~~~~$$
what proves $h$ is strictly trace-convex
along $M$.
The last inequality in the proposition is obtained form (4) that holds for
$\bar{X}=\bar{\nabla} h$ (see theorem 6),
where $\|\bar{X}^{\top}\|=
\|(\bar{\nabla}  h)^{\top} \|\leq \rho\circ F\|\bar{\nabla}\rho\|=\rho\circ F$,
by following \cite{LS}, that we describe now. Since (7) still holds
$\sup_M \rho\circ F$ cannot be achieved, and
given a bounded domain
$D$ of $M$ with boundary $\partial D$ with unit normal $\nu$, we have
$$
A(\partial D)\sup_M (\rho\circ F)\geq \left|\int_{\partial D}
\bar{g}(\bar{X}^{\top},\nu)dA\right| \geq \int_D{\rm div}_g(\bar{X}^{\top})dV
\geq m\alpha V(D).$$
\end{proof}
\noindent
\em Proof of theorem 4. \em  
~This is an immediate consequence of 
previous corollary and theorem 6.\qed

\vspace*{5mm}\noindent
Now we specify for the particular case $\oM=\Sigma'\times \Sigma$,
where $(\Sigma',g_{\Sigma'})$ and $(\Sigma, g_{\Sigma})$ are
Riemannian manifolds of dimension $d'\leq m$ and $d\geq n$ respectively 
where $d+ d'=n+ m$.
Let us fix a point $x_0\in \Sigma'$ 
and denote by $r_{\Sigma'}$
 the distance function in $\Sigma'$ to  $x_0$.
We identify $\Sigma$ with $x_0\times \Sigma$, a totally convex set. 
For 
$(x,y)\in \oM$, we have
$$\rho((x,y))=\bar{d}((x,y), x_0\times \Sigma)=\bar{d}((x,y),(x_0,y))
=d_{\Sigma'}(x,x_0)=
 r_{\Sigma'}(x).$$
Thus,  $ h(x,y)=\ha r^2_{\Sigma'}(x).$
If $F(p)=(x,y)\in \Sigma' \times \Sigma$, and $l=r_{\Sigma'}(x)$ then
$\sigma_{F(p)}(t)=(\sigma^{\Sigma'}(t),y)$ and 
$\sigma'_{F(p)}(l)=((\sigma^{\Sigma'})'(l),0)$
 where
$\sigma^{\Sigma'}$ is a unit geodesic on $\Sigma'$ with $\sigma^{\Sigma'}(0)=
x_0$ and $\sigma^{\Sigma'}(l)=x$.  
Let $\pi_{(x,y)}:T_x\Sigma'=T_y\Sigma^{\bot}\to T_{p}M$, $\pi(v)=v^{\top}$.
 Therefore, 
$F$ is $h$-bounded iff $M$ is immersed in $B_R(x_0)\times \Sigma$
where $B_R(x_0)$ is a ball  in  $\Sigma'$ of radius 
 $R<+\infty$, and
if $\pi$ has sup-norm bounded away from 
zero, then $h$
 is strictly trace-convex on $M$.
\begin{proposition} Let $\Sigma'$ be $m$-dimensional  and $\Sigma$
$n$-dimensional complete connected Riemannian manifolds with
 nonpositive sectional curvatures,  $h:\oM \to [0,+\infty)$,
$h(x,y)=\ha r^2_{\Sigma'}(x)$,  where $r_{\Sigma'}$ is the
distance function in $\Sigma'$ to a given point $x_0$,
and $B_R(x_0)$ the ball of radius $R$ of $\Sigma'$.
If $F:M \to   B_R(x_0)\times \Sigma$ is
a complete minimal 
submanifold  $h$-properly
 immersed and  there exist a constant $C>0$ such that
$M$ is locally
 the graph of a local map $f:B_R(x_0)\to \Sigma$ with $f^*g_{\Sigma}\leq
Cg_{\Sigma'}$,
then $M$
has pure point spectrum.
\end{proposition}
\begin{proof}  First we note that 
the sectional curvature of $\oM$ is also nonnegative.
In the  particular case $d'=m$, and if locally
$M$ is  the graph of a local  map $f:\Sigma'\to \Sigma$,
then we show that the trace-convexity holds if 
$f^*g_{\Sigma}\leq C g_{\Sigma'}$, 
for some constant $C>0$. At a given pointi $p\in M$, let 
$\lambda_1^2\geq \ldots\geq \lambda_m^2$ be the eigenvalues
of $f^*g_{\Sigma}$ with corresponding $g_{\Sigma'}$-o.n.\ basis $a_i$
of eigenvectors.
Then it follows that
at $F(p)=(x,y)=(x,f(x))$,   $df_x(a_i)=\lambda_i a_{i+m}$, 
where $a_i,a_{\alpha}$, $i=1,\ldots, m$, 
$\alpha=m+1\ldots,m+n$ defines an o.n.\ basis of $T_{(x,y)}\oM$
(note that $\lambda_i=0$ for $i>\min \{m,n\}$, so we can always find
such basis).
Then $e_i=(a_i +\lambda_i a_{i+m})/(1+ \lambda_i^2)^{1/2}$
constitutes an o.n.\ basis of the graph and 
$$\|(\sigma'_{F(p)}(l))^{\top}\|^2= \|((\sigma^{\Sigma'})'(l),0)^{\top}\|^2=
\sum_i |g_{\Sigma'}((\sigma^{\Sigma'})'(l),a_i)|^2/{(1+\lambda_i^2)}\geq 
 \frac{1}{1+C},$$
and the proposition is proved. 
\end{proof}

\noindent
\em Remark. \em 
The previous proposition should be compared with a similar result
for the case $\Sigma$ and $\Sigma'$ Euclidean spaces in \cite{BJM}.
If $\Sigma'=\mathbb{R}^m$ and $\Sigma=\mathbb{R}^{n}$,
according to \cite{ABD}
the immersion in the previous proposition cannot be
properly immersed in $\mathbb{R}^{m+n}$, if $m\geq n+1$.


\begin{thebibliography}{99}
{
\bibitem{ABD} L.J.\ Al\'{\i}as,  G.\ Pacelli Bessa,  M.\ Dajczer, 
  \em Counterexamples to Calabi conjectures on minimal
         hypersurfaces cannot be proper, 
\em arXiv/math:0812.0623v1
\bibitem{BM} G.P.\ Bessa  and J.F. Montenegro, 
\em An extension of Barta's theorem and geometric applications.  \em
Ann.\ Global Anal.\ Geom.\ {\bf 31} (2007), no.4, 345-362.
\bibitem{BJM}   G.P.\ Bessa, L.P. Jorge and J.F. Montenegro,
\em On the spectrum of the Martin-Morales minimal surfaces. \em
arXiv/math:0809.1173.
\bibitem{Ca} E. Calabi,  \em Problems in differential geometry, \em 
Ed. S. Kobayashi and J. Eells, Jr., Proceedings of the
United States-Japan Seminar in Differential Geometry, Kyoto, 
Japan, 1965. Nippon Hyoronsha Co., Ltd., Tokyo (1966) 170.
\bibitem{Che} J.\ Cheeger, \em A lower bound for the smallest eigenvalue
of the Laplacian. \em  "Problems in Analysis", 195-199. Princeton Univ.\
Press, Princeton, New Jersey, 1970.
\bibitem{CM} T.H.\ Colding and W.P. Minicozzi II, \em The Calabi-Yau conjectures
for embedded surfaces, \em  Annals of Math. {\bf 161} (2005) 727-758.
\bibitem{Donn}  H.\ Donnelly, \em On the essential spectrum of a complete
Riemannian manifold. \em Topology {\bf 20} (1981), 1-14. 
\bibitem{DL} H.\ Donnelly and P.\ Li, \em
Pure point spectrum and negative curvature for noncompact manifolds. \em
 Duke Math.\ J.\ {\bf 46} (1979), 497–503.
\bibitem{Gray} A.\ Gray, \em Tubes.\em  Second edition. 
Progress in Mathematics, {\bf 221} (2004), Birkh\"{a}user Verlag, Basel. 
\bibitem{Ka} A.\ Kasue, \em On Laplacian and hessian comparison theorems.
\em Proc.\ Japan Acad.\ {\bf 58} Ser. A (1982), 25-28.
\bibitem{LS} G.\ Li and I. Salavessa, \em Bernstein-Heinz-Chern results
in calibrated manifolds. \em  arXiv:0802.0946. 
\bibitem{MM1} F.\ Mart\'{\i}n ́ and S.\ Morales, 
\em Complete proper minimal surfaces in convex bodies of $R^3$. \em 
Duke Math.\ J.\  {\bf 128}, (2005), 559–593.
\bibitem{MM2} F.\ Mart\'{\i}n ́ and S.\ Morales, 
\em Complete proper minimal surfaces in convex bodies of $R^3$. II. 
The behavior of the limit set. \em 
Comment.\ Math.\ Helv.\ 81, (2006), 699-725.
\bibitem{Na} N.\ Nadirashvili, \em Hadamard’s and Calabi-Yau’s 
conjectures on negatively curved and minimal surfaces. \em 
Invent.\ Math., {\bf 126}, (1996), 457–465.
\bibitem{R-M} A.\ Ranjbar-Motlagh, \em Generalizations of the Liouville
theorem. \em Diff.\ Geom.\ Appl.\ {\bf 26} (2008), 339-345.
\bibitem{RS} M.\ Rigoli and I.\ Salavessa, 
\em Conformal and Isometric immersions of Riemannian
manifolds.\em  Math.\ Z {\bf 196} (1987), 293-300.
\bibitem{S0} I.M.C.\ Salavessa, \em Graphs with parallel mean curvature 
and a variational problem  in conformal geometry. \em 
Ph.D. Thesis, University of Warwick, 1987
}

\end{thebibliography}
\end{document}